\magnification\magstep1
\hsize 15.5truecm
\scrollmode

 at 17.28truept
 at 14.4truept
 at 12truept
 at 10.95truept
 at 17.28truept
 at 10truept

\def\title#1{%
\vskip0pt plus.3\vsize\penalty-100%
\vskip0pt plus-.3\vsize\bigskip\vskip\parskip%
\bigbreak\bigbreak\centerline{\bf #1}\bigskip%
}

\def\chapter#1#2{\vfill\eject
\centerline{\bf Chapter #1}
\vskip 6truept%
\centerline{\bf #2}%
\vskip 2 true cm}

\def\section#1#2{%
\def\\{#2}%
\vskip0pt plus.3\vsize\penalty-100%
\vskip0pt plus-.3\vsize\bigskip\vskip\parskip%
\par\noindent{\bf #1\hskip 6truept%
\ifx\empty\\{\relax}\else{\bf #2\smallskip}\fi}}

\def\subsection#1#2{%
\def\\{#2}%
\vskip0pt plus.3\vsize\penalty-20%
\vskip0pt plus-.3\vsize\medskip\vskip\parskip%
\def\TEST{#1}%
\noindent{\ifx\TEST\empty\relax\else\bf #1\hskip 6truept\fi%
\ifx\empty\\{\relax}\else{#2\smallskip}\fi}}

\def\proclaim#1{\medbreak\begingroup\noindent{\bf #1.---}\enspace\sl}

\def\endproclaim{\endgroup\par\medbreak}

\def\demo{\smallbreak\begingroup\noindent{\bf Proof.\ }\rm}

\def\lbreak{\hfil\break}


\def\comfig#1#2\par{
\medskip
\centerline{\hbox{\hsize=10cm\eightpoint\baselineskip=10pt
\vbox{\noindent #1}}}\par\centerline{ Figure #2}}

\def\figcom#1#2\par{
\medskip
\centerline
{Figure #1}
\par\centerline{\hbox{\hsize=10cm\eightpoint\baselineskip=10pt
\vbox{\noindent #2}}}}
\def\bull{~\vrule height .9ex width .8ex depth -.1ex}


\def\lbreak{\hfill\break}

\def\comfig#1#2\par{
\medskip
\centerline{\hbox{\hsize=10cm\eightpoint\baselineskip=10pt
\vbox{\noindent{\sl  #1}}}}\par\centerline{{\bf Figure #2}}}

\def\figcom#1#2\par{
\medskip
\centerline
{{\bf Figure #1}}
\par\centerline{\hbox{\hsize=10cm\eightpoint\baselineskip=10pt
\vbox{\noindent{\sl  #2}}}}}

\def\em{\sl}

\def\\{\hfill\break}

\def\lbreak{\hfill\break}

\def\bull{~\vrule height .9ex width .8ex depth -.1ex}









\def\diagram#1{%
\def\normalbaselines{\baselineskip20pt\lineskip3pt\lineskiplimit3pt}%
\matrix{#1}}

\def\a{\alpha}
\def\b{\beta}

\def\la{\lambda}

\def\CC{{\bf C}}

\def\PP{{\bf P}}

\def\ZZ{{\bf Z}}

\def\cD{{\cal D}}
\def\cE{{\cal E}}
\def\cF{{\cal F}}

\def\cL{{\cal L}}
\def\cM{{\cal M}}

\def\cO{{\cal O}}

\def\cZ{{\cal Z}}

\def\Hom{\mathop{\rm Hom}\limits}


\def\la{\lambda}

\def\sqr#1#2{{\vcenter{\hrule height.#2pt%
\hbox{\vrule width.#2pt height#1pt\kern#1pt%
\vrule width.#2pt}%
\hrule height.#2pt}}}

\newfam\gotfam
\font\twlgot=eufm10 at 12pt
\font\tengot=eufm10

\font\sevengot=eufm7
\textfont\gotfam=\twlgot
\scriptfont\gotfam=\tengot 
\scriptscriptfont\gotfam=\sevengot


\input amssym.def
\input amssym.tex

\overfullrule=0pt
\input psfig
\input psbox
\psfordvips
\def\pr{\partial}
\def\diag{{\rm diag}}\def\Ker{{\rm Ker}}
\def\refAS{[1]}
\def\refBD{[2]}
\def\refD{[3]}
\def\refER{[4]}
\def\refFG{[5]}
\def\refFFR{[6]}
\def\refFV{[7]}
\def\refF{[8]}
\def\refN{[9]}
\def\refS{[10]}

\centerline{\bf Separation of variables for Gaudin-Calogero systems}
\bigskip
\centerline{B. Enriquez, B. Feigin and V. Rubtsov}
\medskip
{\bf Abstract.} {\em
We construct an elliptic analogue of Sklyanin's separation of variables
for the $sl(2)$ Gaudin system, using an adaptation of Drinfeld's Radon
transformations. 
}
\medskip

{\bf Introduction.} 

The geometric Langlands conjectures, as formulated  by Beilinson and
Drinfeld in \refBD, aim at the construction of certain $\cD$-modules on
the moduli space of $G$-bundles over (punctured) curves ($G$ a reductive
group). Positive characteristic versions of these conjectures had been
solved earlier by Drinfeld in \refD, in the case $G=GL_{2}$. 

In the paper \refFFR, the $\cD$-modules arising from the construction of
\refBD\ were studied in the special case of a rational curve with marked
points, and identified with the Gaudin model. Then,
in \refF, Drinfeld's construction of local systems on the
moduli space of rank two vector bundles on a curve in positive
characteritic (\refD) was adapted to the complex situation. The identification
of these two constructions amounts to Sklyanin's separation of variables
(\refS), as it was noticed in \refF. This computation is recalled in the first
part of this text. 

The question has been raised in \refF\ to construct a similar separation of
variables for the Gaudin-Calogero systems, which were computed in \refER\ and
\refN, and play a similar role
in the case of a punctured elliptic curve. This note aims at solving
this question. In the present case Drinfeld's diagrams for Radon
transformation have to be slightly modified. 

It is also worth to note that the systems presented here, are the
specialization at the critical level, of the
Knizhnik-Zamolodchikov-Bernard equations on the torus. According to
the general viewpoint that the diagrams of \refD\ are related to the
Drinfeld-Sokolov reduction, the generalisation
of the present work to noncritical level should relate these equations
to the Virasoro correlators on the torus. In \refFG, the Bethe equations
were connected with the unitarity property of the KZB equation on the
torus.
It would be desirable to clarify further the connection between these
issues.  

We express our thanks to A. Stoyanovsky and A. Varchenko for having
discussed with us the
content of this paper, and to C. Sabbah for consultations on
$\cD$-modules. V.R. was supported by the INTAS grants 93-2494 and
1010-CT93-0023, and by the CNRS; he expresses his thanks to these
institutions.  

\section{1.}{Separation of variables for $sl_{2}$ Gaudin systems.}

Let us begin with some reminders on the Gaudin system. Let $X=\CC
P^{1}$, $z_{\a}$ be marked points on $X$, $\a=1,\cdots,N$, 
$G$ be $SL_{2}(\CC)$, $B\subset G$ be the upper
triangular subgroup. The moduli space $\cM_{G}(X, z_{\a})$ of
$G$-bundles on $\CC P^{1}$, with parabolic structures at $z_{\a}$, is the
disjoint union of the $\cM_{G}^{(n)}(X, z_{\a})$, $n\ge 0$,
corresponding to the parabolic structures on the sheaf
$\cO(n\infty)\oplus\cO(-n\infty)$. We then identify  $\cM_{G}^{(n)}(X,
z_{\a})$ with $P_{n}\backslash (G/B)^{N}$, $P_{0}=G$, and $P_{n}=\{
\pmatrix{ t& p(z)\cr 0 & t^{-1}\cr }, t\in\CC^{\times}, p(z)\in\CC[z],
{\rm deg}(p)\le 2n\}$ for $n>0$, $P_{0}$ acting diagonally and $P_{n}$ 
acting on
the $\a$-th factor by left translation, after the replacement of $z$ by
$z_{\a}$. In what follows we will deal with $\cM_{G}^{(0)}(X, z_{\a})$.

Let us fix weights, $\la_{\a}$, $\a=1,\cdots,N$. On $\cM_{G}^{(0)}(X,
z_{\a})$ lives the bundle $\cL_{(\la_{\a})}$, quotient of the bundle
$\boxtimes_{\a=1}^{N}\cL_{\la_{\a}}$ ($\cL_{\la_{\a}}$ is the line
bundle on
$G/B$, corresponding to the weight $\la_{\a}$). 
The natural action of $Z(U_{-2}\hat{sl}_{2})_{loc}$ (the center of the
local completion of the enveloping algebra of the central extension of
$sl_{2}(\CC((z)))$ at level $-2$) is by differential operators $L_{\a}$,
$\a=1,\cdots,N$, which were identified in \refFFR\ with the Gaudin
hamiltonians 
$L_{\a}=\sum_{\b\ne\a}{{I_{a}^{(\a)}I^{a(\b)}}\over{z_{\a}-z_{\b}}}$,
$I_{a}$, $I^{a}$ being an orthogonal basis of $sl_{2}(\CC)$. 

Following the conjectures of \refBD, the $\cD$-modules on 
$\cM_{G}(X, z_{\a})$ (twisted by $\cL_{(\la_{\a})}$) defined by 
$L_{\a}-\mu_{\a}$ should satisfy the Hecke eigenvalue property. In \refD, a
construction of such modules was given in the case of a curve of genus
$>0$ without punctures. 

Following \refF, let us show how Drinfeld's construction
in \refD\ can be adapted for $X= \CC P^{1}$ with marked points $z_{\a}$. 
Consider the space 
$$
\eqalign{
\cM^{(0)}_{B}(X,z_{\a}) = \{({\rm parabolic\ }&{\rm \ structure\ on\ }
\cO_{\CC P^{1}}^{2} {\rm \  at\ } z_{\a}, \cr & 
{\rm class\  of\  morphisms\ }
\cO_{\CC P^{1}}\to \cO_{\CC P^{1}}^{2})\},
}
$$
 the morphisms being considered
up to $\cO_{\CC P^{1}}$-automorphisms; the cokernel of the morphism
considered in this definition is $\cO_{\CC P^{1}}$. It is natural to
consider $K_{\CC P^{1},z_{\a}} 
= K_{\CC P^{1}}(\sum_{\a=1}^{N}(z_{\a}))$ as the
canonical bundle in our punctured situation, and then the space 
$\{$classes of morphisms $K^{-1}_{\CC P^{1},z_{\a}}
\to \cO_{\CC P^{1}}\}$, its mapping $\pi$ to $X^{N-2}$ (given by the
zeroes of a given section) and the diagram
$$
\diagram{
&&& \cZ &&&\cr
&& p_{{0}}\swarrow && \searrow q_{{0}}&&\cr
& \cM_{B}^{(0)}(\CC P^{1},z_{\a})&&& \PP\Hom (K_{\CC P^{1},z_{\a}}^{-1},
\cO_{\CC P^{1}})& \buildrel\pi\over\longrightarrow X^{(N-2)}&	 \cr
  p \swarrow &&\searrow && \swarrow  &&\cr
\cM^{(0)}_{G}(\CC P^{1},z_{\a}) &&& \{ \cO_{\CC P^{1}}\}&&&\cr}\leqno(1)
$$
$p$ being the projection on the first factor, and the correspondance
$\cZ$ being defined to be the set of $((l_{\a},i),j)$, $l_{\a}$: line
in the fiber of $\cO_{\CC P^{1}}^{2}$ at $z_{\a}$, $i$: morphism
$\cO_{\CC P^{1}} \to \cO_{\CC P^{1}}^{2}$, $j$: morphism $K_{\CC
P^{1},z_{\a}}\to\cO_{\CC P^{1}}$, proportional $i,j$'s
being considered equivalent, such that denoting by $k:\cO_{\CC
P^{1}}^{2}\to \cO_{\CC P^{1}}$, the cokernel mapping of $i$, there exists
a lift $j': K_{\CC P^{1},z_{\a}}\to \cO_{\CC
P^{1}}^{2}$ of $j$ (i.e., we have $j=k\circ j$), compatible with the
parabolic structure (i.e., the image of $j'$ at $z_{\a}$ should be the
line $l_{\a}$).

Let us fix weights $\lambda_{\a}$ for each $\a$; the $\cD$-modules we
will consider will be twisted by the following line bundles:
$\cL_{(\lambda_{\a})}$ on $\cM_{G}^{(0)}(X,z_{\a})$,
$p^{*}\cL_{(\lambda_{\a})}$ on $\cM_{B}^{(0)}(X,z_{\a})$, and
$p_{0}^{*}p^{*}\cL_{(\lambda_{\a})}$ on $\cZ$. For $Y$ a variety and
$\cL$ a line bundle on $Y$, we denote $(\cD_{Y})_{\cL}=
\cL\otimes\cD_{Y}\otimes\cL^{-1}$. 
Let us fix now complex numbers $\mu_{\a}, \a=1,\cdots,N$, s.t.
$\sum_{\a=1}^{N}\mu_{\a}=0$,
$\sum_{\a=1}^{N}\mu_{\a}z_{\a}+\sum_{\a=1}^{N}2\la_{\a}(\la_{\a}-1)=0$, 
$\sum_{\a=1}^{N}\mu_{\a}z_{\a}^{2}+\sum_{\a=1}^{N}4\la_{\a}(\la_{\a}-1)
z_{\a}=0$; we associate to
them the operator on $X$, 
$$
D_{(\lambda_{\a}),(\mu_{\a})}=
2\partial_{w}^{2}-\sum_{\a=1}^{N}{\mu_{\a}\over w-z_{\a}} -
\sum_{\a=1}^{N}{2\lambda_{\a}(\lambda_{\a}-1)\over(w-z_{\a})^{2}},
\leqno(2)
$$
and the $\cD_{X}$-module $\cE_{(\lambda_{\a}),(\mu_{\a})}= \cD_{X}/
\cD_{X}D_{(\lambda_{\a}),(\mu_{\a})}$.

Consider on the other hand on $\cM_{G}^{(0)}(X,z_{\a})$ the twisted
$\cD$-module 
$$
\cM_{(\mu_{\a})}=
\Big(\cD_{\cM_{G}^{(0)}(X,z_{\a})}\Big)_{\cL_{(\la_{\a})}}
/  \sum_{\a=1}^{N}\Big(\cD_{\cM_{G}^{(0)}(X,z_{\a})}
\Big)_{\cL_{(\la_{\a})}}
(L_{\a}-\mu_{\a}) 
$$ 
(the conditions on $\mu_{\a}$ correspond to the
relations on the $L_{\a}$, 
$$
\sum_{\a=1}^{N}L_{\a}=0,
\sum_{\a=1}^{N}L_{\a}z_{\a}+\sum_{\a=1}^{N}2\la_{\a}(\la_{\a}-1)
=ef+fe+{1\over 2}h^{2},
$$ 
$$
\sum_{\a=1}^{N}L_{\a}z_{\a}^{2}+\sum_{\a=1}^{N}4\la_{\a}(\la_{\a}-1)
z_{\a}=2(e_{1}f+f_{1}e+{1\over 2}h_{1}h)),
$$ 
$e=\sum_{\a=1}^{N}e^{(\a)}$, 
$e_{1}=\sum_{\a=1}^{N}z_{\a}e^{(\a)}$, analogous relations for $f$,
$f_{1}$, $h$, $h_{1}$). We would like to show: 

\proclaim{Proposition} (cf. \refF.)
There is a homomorphism of $\cD$-modules
$$
\pi^{*}\cE^{(N-2)}_{(\la_{\a}),(\mu_{\a})}\to
R(q_{0})_{*}p_{1}^{*}p^{*}
\cM_{(\mu_{\a})}[N],
$$ 
which is an isomorphism over $\pi^{-1}((X-\{
\infty\})^{(N-2)}-\Delta)$
($\Delta$ is the diagonal part of
$(X-\{\infty\})^{(N-2)}$). \endproclaim 

(Here we denote, for $\cF$ a sheaf on a manifold $V$, by $\cF^{(n)}$ the
sheaf $(p_{V})_{*}(\cF^{\boxtimes n})$ on $V^{(n)}=V^{n}/S_{n}$, $p_{V}$
being the projection $V^{n}\to V^{(n)}$.)

\demo Let us give coordinates to the spaces of diagram (1).  
$$
\cM_{G}^{(0)}(X,z_{\a})\simeq G\backslash (G/B)^{N}=G\backslash (\CC
P^{1})^{N},
$$ 
choosing the identification $\CC P^{1}\simeq G/B$,
$t\mapsto\pmatrix{ 1 & 0\cr t^{-1} & 1 \cr} B$, 
$0\mapsto\pmatrix{ 0 & -1\cr 1 & 0 \cr} B$; then $G$  acts on $(\CC
P^{1})^{N}$ by homographic transformations. Now
$\cM_{B}^{(0)}(X,z_{\a})=B\backslash (G/B)^{N}$; after fixing $\cO_{\CC
P^{1}}\to \cO_{\CC P^{1}}^{2}$ to be $(1,0)$, the lines $l_{\a}$ are
$\CC(1,t_{\a}^{-1})$, the $t_{\a}^{-1}$ being defined up to a global
affine transformation. An element of $\Hom(K^{-1}_{\CC P^{1},z_{\a}}, 
\cO_{\CC P^{1}})$ is a $1$-form
$\sum_{\a=1}^{N}{u_{\a}dz\over {z-z_{\a}}}
$, with $\sum_{\a=1}^{N}u_{\a}=0$.
The incidence relation defining $\cZ$ is
$\sum_{\a=1}^{N}u_{\a}t_{\a}=0$, since the first component of $j'$ has
to be $\sum_{\a=1}^{N}{u_{\a}t_{\a}dz\over {z-z_{\a}}}$, and should be
regular at $\infty$. The map $\pi$ associates to $(u_{\a})$, the
solutions $(w_{i})$ of $\sum_{\a=1}^{N}{u_{\a}\over{z-z_{\a}}}=0$
(counting $k$ times $\infty$, if this function is $\sim c/w^{2+k}$ for
$w\to \infty$,
$c\ne 0$). 

Let $p_{1}$ be the natural projection of $(G/B)^{N}$ on $G\backslash
(G/B)^{N}$, then
$$
p_{1}^{*}(\cM_{(\mu_{\a})})=\cD_{(G/B)^{N}}/
\sum_{\a}\cD_{(G/B)^{N}}(L_{\a}-\mu_{\a})+\cD_{(G/B)^{N}}sl_{2}(\CC).
$$

Introduce the formal variable $z$, then 
$$
\sum_{\a=1}^{N}{{L_{\a}-\mu_{\a}}\over{z-z_{\a}}}+
\sum_{\a=1}^{N}{{2\la_{\a}(\la_{\a}-1)}\over{(z-z_{\a})^{2}}}
=
e(z)f(z)+f(z)e(z)+{1\over 2}h(z)^{2}-
\sum_{\a=1}^{N}{\mu_{\a}\over{z-z_{\a}}}
$$
with $e(z)=\sum_{\a=1}^{N}{e^{(\a)}\over {z-z_{\a}}}$, etc., 
$e^{(\a)}=t_{\a}^{2}{\pr\over{\pr t_{\a}}}
+2\la_{\a}t_{\a}$, $f^{(\a)}=-{\pr\over{\pr t_{\a}}}$, 
$h^{(\a)}=2(t_{\a}{\pr\over{\pr t_{\a}}}+\la_{\a})$.
The Radon transform
of the $\cD$-module generated by the $L_{\a}-\mu_{\a}$ is the
$\cD$-module generated by the $\bar{L}_{\a}-\mu_{\a}$, where
$$
\sum_{\a=1}^{N}{{\bar{L}_{\a}-\mu_{\a}}\over{z-z_{\a}}}+
\sum_{\a=1}^{N}{{2\la_{\a}(\la_{\a}-1)}\over{(z-z_{\a})^{2}}} =
\bar{e}(z)\bar{f}(z)+\bar{f}(z)\bar{e}(z)+{1\over 2}\bar{h}(z)^{2}-
\sum_{\a=1}^{N}{\mu_{\a}\over{z-z_{\a}}}
$$
$\bar{e}(z)=\sum_{\a=1}^{N}{\bar{e}^{(\a)}\over {z-z_{\a}}}$, analogous
formulae for $\bar f(z)$, $\bar h(z)$,  
$
\bar{e}^{(\a)}= - ( u_{\a}({\pr\over{\pr u_{\a}}})^{2}
+2(\la_{\a}+1){\pr\over {\pr u_{\a}}} ) , 
\bar{f}^{(\a)}=u_{\a},$ $ 
\bar{h}^{(\a)}= - 2(u_{\a}{\pr\over{\pr u_{\a}}}+\la_{\a}+1).$

Consider the operator
$$
\hat{L}(w_{i})=\sum_{\a=1}^{N}{1\over{w_{i}-z_{\a}}}(\bar{L}_{\a}-\mu_{\a}),
$$
and let
$\hat{e}(w_{i})=\sum_{\a=1}^{N}{1\over{w_{i}-z_{\a}}}\bar{e}^{(\a)}$,
analogous formulae for $\hat{f}(w_{i})$, $\hat{h}(w_{i})$. Then 
$$
\eqalign{
\hat{L}(w_{i}) & +\sum_{\a=1}^{N}{\mu_{\a}\over{w_{i}-z_{\a}}}  
+\sum_{\a=1}^{N}{{2\la_{\a}(\la_{\a}-1)}\over{(w_{i}-z_{\a})^{2}}}
-[\hat{e}(w_{i})\hat{f}(w_{i})+\hat{f}(w_{i})
\hat{e}(w_{i})+{1\over 2}\hat{h}(w_{i})^{2}] \cr &
=-\sum_{1\le\a,\b\le N}
{1\over{w_{i}-z_{\a}}}\{
[\bar{e}^{(\a)},{1\over {w_{i}-z_{\b}}}]\bar{f}^{(\b)}
+[\bar{f}^{(\a)},{1\over {w_{i}-z_{\b}}}]\bar{e}^{(\b)}
\cr &
+{1\over 2}[\bar{h}^{(\a)},{1\over {w_{i}-z_{\b}}}]\bar{h}^{(\b)}\} \cr
&=
-\sum_{1\le\a,\b\le N} - {1\over{w_{i}-z_{\a}}}[u_{\a}{\pr\over{\pr
u_{\a}}}
,{\pr\over{\pr
u_{\a}}}({1\over {w_{i}-z_{\b}}})]_{+}\cdot u_{\b}
\cr &
+{2\over{w_{i}-z_{\a}}}u_{\a}{\pr\over{\pr u_{\a}}}({1\over
{w_{i}-z_{\b}}})(u_{\b}{\pr\over{\pr u_{\b}}}) 
\cr
& -\sum_{1\le \a,\b \le N} - {1\over{w_{i}-z_{\a}}}
\{[2(\la_{\a}+1){\pr\over{\pr u_{\a}}}
, {1\over{w_{i}-z_{\b}}}]u_{\b}
\cr &
-{2\over{w_{i}-z_{\a}}}u_{\a}{\pr\over{\pr
u_{\a}}}({1\over{w_{i}-z_{\b}}})(\la_{\b}+1)\}
\cr}
$$
with $[a,b]_{+}=ab+ba$. 
Now,
$$
\sum_{\b=1}^{N} {\pr\over{\pr
u_{\a}}}\Big({1\over{w_{i}-z_{\b}}}\Big)u_{\b}
=-{1\over {w_{i}-z_{\a}}},
$$
$$
\sum_{\b=1}^{N} u_{\a}\Big({\pr\over{\pr
u_{\a}}}\Big)^{2}\Big({1\over{w_{i}-z_{\b}}}\Big)u_{\b}
=-2u_{\a}{\pr\over{\pr u_{\a}}}\Big({1\over {w_{i}-z_{\a}}}\Big),
$$
so the last line gives zero, and the term in $[\  ,\  ]_{+}$ gives
$$
\eqalign{
-2\sum_{\a=1}^{N}{1\over{w_{i}-z_{\a}}}u_{\a}{\pr\over{\pr
u_{\a}}}\cdot{1\over{w_{i}-u_{\a}}}
& +2\sum_{\a=1}^{N}{1\over{w_{i}-z_{\a}}}u_{\a}{\pr\over{\pr
u_{\a}}}\Big({1\over{w_{i}-u_{\a}}}\Big)
 \cr &=
-2\sum_{\a=1}^{N}{1\over{(w_{i}-z_{\a})^{2}}}u_{\a}{\pr\over{\pr
u_{\a}}}.
}
$$
(the dot denotes the product of differential operators). 
Then we deduce from
$$
C{{\prod_{i=1}^{N-1}(z-w_{i})}\over{\prod_{\a=1}^{N}(z-z_{\a})}}= 
\sum_{\a=1}^{N}{u_{\a}\over{z-z_{\a}}},
$$
$({dC\over C}+\sum_{i=1}^{N-2}{{dw_{i}}\over{w_{i}-z}})
\sum_{\a=1}^{N}{u_{\a}\over{z-z_{\a}}}=
\sum_{\a=1}^{N}{du_{\a}\over{z-z_{\a}}}$, so 
$du_{\a}=u_{\a}({dC\over C}+\sum_{i=1}^{N-2}{{dw_{i}}\over{w_{i}-z_{\a}}})$
and 
$$
{\pr\over{\pr
w_{i}}}=\sum_{\a=1}^{N}{u_{\a}\over{w_{i}-z_{\a}}}{\pr\over{\pr
u_{\a}}};
$$
so the remaining term gives 
$-2\sum_{\b=1}^{N}{\pr\over{\pr
w_{i}}}({1\over{w_{i}-z_{\b}}})(u_{\b}{\pr\over{\pr u_{\b}}})
=
\sum_{\b=1}^{N}{2\over{(w_{i}-z_{\b})^{2}}}
(u_{\b}{\pr\over{\pr u_{\b}}})$. 

Finally, 
$$\eqalign{
\hat{L}(w_{i})=
\hat{e}(w_{i})\hat{f}(w_{i})+\hat{f}(w_{i})
\hat{e}(w_{i})+{1\over 2}\hat{h}(w_{i})^{2} &
-\sum_{\a=1}^{N}{\mu_{\a}\over{w_{i}-z_{\a}}}
-\sum_{\a=1}^{N}{ {2\la_{\a}(\la_{\a}-1)} \over{(w_{i}-z_{\a})^{2}}}.\cr}
$$
Now, $\hat{f}(w_{i})=0$, and $\hat{h}(w_{i})=-2[{\pr\over{\pr
w_{i}}}+A(w_{i)}]$, with
$A(w_{i})=\sum_{\a=1}^{N}{{\la_{\a}+1}\over{w_{i}-z_{\a}}}$, so 
$$
\hat{L}(w_{i})=2\Big({\pr\over{\pr w_{i}}}+A(w_{i})\Big)^{2}
-\sum_{\a=1}^{N}{\mu_{\a}
\over{w_{i}-z_{\a}}}
-\sum_{\a=1}^{N}{ {2\la_{\a}(\la_{\a}-1)} \over{(w_{i}-z_{\a})^{2}}}.
$$
In this way, we have constructed an epimorphism from the $(N-2)$-th symmetric
power of $\cE_{(\la_{\a}),(\mu_{\a})}$ (restricted to the complement of
diagonals) to the $\cD$-module generated by the $\bar L_{\a}-\mu_{\a}$'s
(restricted to the complement of the discriminant), and so to the
$\cD$-module generated by the $\bar L_{\a}-\mu_{\a}$'s and the action
of $sl_{2}(\CC)$. 

Let us now show that it induces an isomorphism of the sheaves of local
analytic solutions. Let us start with a
local homomorphism of the first sheaf to the (analytic) 
structure sheaf. It is some
function $(\psi(w_{i}))_{1\le i\le N-2}$, such that
$(\pr_{w_{j}}^{2}+q(w_{j}))\psi(w_{i})=0$, for all $j$. We deduce from
that relation, using our previous computations, 
$$
\sum_{i=1}^{N}{{(\bar L_{\a}-\mu_{\a})\psi}\over{w_{i}-z_{\a}}}=0, 
\quad i=1,\cdots,N-2.\leqno{(3)}
$$
We will consider $\psi$ as a distribution on the space of all
$(u_{i})_{1\le i\le N}$, supported on the hyperplane
$\sum_{i=1}^{N}u_{i}=0$, and analytic on this hyperplane. 
We obtain from (3) 
$$
(\bar L_{\a}-\mu_{\a})\psi=u_{\a}\phi+z_{\a}u_{\a}\rho,\leqno{(4)}
$$
$\phi$ and $\rho$ being distributions of the same nature as $\psi$;
indeed, $\phi$ and $\rho$ can be obtained solving a Cramer system (since
everywhere on the support of $\psi$, we can find two indices $\a\ne\b$
such that $u_{\a}u_{\b}\ne 0$), and 
$\bar L_{\a}$'s commute with $\sum_{\a=1}^{N}u_{\a}$. 
From the relation $\sum_{\a=1}^{N}\bar L_{\a}=0$ follows that $\rho=0$;
from  
$\sum_{\a=1}^{N}\bar L_{\a}z_{\a}+\sum_{\a=1}^{N}2\la_{\a}(\la_{\a}-1)
=\bar e\bar f+\bar f\bar e+{1\over 2}\bar h^{2}$ and 
$\sum_{\a=1}^{N}\bar L_{\a}z_{\a}^{2}+\sum_{\a=1}^{N}4\la_{\a}(\la_{\a}-1)
z_{\a}=2(\bar e_{1}\bar f+\bar f_{1}\bar e+{1\over 2}\bar h_{1}\bar
h)$, follows that  $(\sum_{\a=1}^{N}u_{\a})\bar
e\psi=(\sum_{\a=1}^{N}z_{\a}u_{\a})\phi$, and 
$(\sum_{\a=1}^{N}z_{\a}u_{\a})\bar
e\psi=(\sum_{\a=1}^{N}z_{\a}^{2}u_{\a})\phi$. 
So, on the complement of $\{(u_{\a})|\sum_{\a=1}^{N}z_{\a}u_{\a}=0$ or $
\sum_{\a=1}^{N}z_{\a}^{2}u_{\a}=0$ or
$(\sum_{\a=1}^{N}z_{\a}u_{\a})^{2}-(\sum_{\a=1}^{N}u_{\a})
(\sum_{\a=1}^{N}z_{\a}^{2}u_{\a})\}$, $\phi$ and $\bar e \psi$ will
vanish. So these distributions would have to be supported on a
subvariety of the set of all $(u_{\a})$ of codimension $\ge 2$, which is
impossible.  So we will have 
$$
(\bar L_{\a}-\mu_{\a})\psi=0, \quad \bar e \psi=0.\leqno{(5)}
$$
Since $\bar f\psi=\bar h\psi=0$ by construction, we have shown that
$\psi$ can be considered as a local homomorphism of 
the $\cD$-module generated by the $\bar L_{\a}-\mu_{\a}$
(restricted to the complement of the discriminant) to the structure
sheaf. 

These two morphisms are clearly inverse to each other; if we show that both
$\cD$-modules have their characteristic varieties supported on the zero
section, this will prove the proposition.  

For $\cE^{(N-2)}_{(\la_{\a}),(\mu_{\a})}$ it is clear, since it is true
for $\cE_{(\la_{\a}),(\mu_{\a})}$. On the other hand, the characteristic
variety of $R(q_{0})_{*}p_{1}^{*}p^{*}
\cM_{(\mu_{\a})}[N]$ is the set of $(u_{i}, \xi_{i})$, with
$\sum_{i=1}^{N}u_{i}=0$ and up to equivalence $(u_{i},\xi_{i})
\sim (u_{i},\xi_{i}+\la)$, and 
$$
\sum_{i=1}^{N}{u_{i}\over{z-z_{i}}}=
\Big(\sum_{i=1}^{N}{{u_{i}\xi_{i}}\over{z-z_{i}}}\Big)
\Big(\sum_{i=1}^{N}{{u_{i}\xi_{i}^{2}}\over{z-z_{i}}}\Big).
$$
This equation gives 
$$
\sum_{i=1}^{N}{{u_{i}}\over{z-z_{i}}}
={{RA^{2}(z)}\over{\prod_{i=1}^{N}(z-z_{i})}}, 
\sum_{i=1}^{N}{{u_{i}\xi_{i}}\over{z-z_{i}}}
={{RAB(z)}\over{\prod_{i=1}^{N}(z-z_{i})}},
\sum_{i=1}^{N}{{u_{i}\xi_{i}^{2}}\over{z-z_{i}}}
={{RB^{2}(z)}\over{\prod_{i=1}^{N}(z-z_{i})}},
$$
$R$, $A$, $B$
polynomials. The set $\{w_{i}, 1\le i\le N-2\}$ is the union of the set
of zeroes of $RA^{2}$, and of 
$\infty$ counted $N-2-{\rm deg} RA^{2}$ times. Since the $w_{i}$'s are
pairwise distinct, we have $A=$const., and deg$R=N-2$ or $N-3$. Since no
$w_{i}$ coincides with $\infty$, deg$R=N-2$; so, $B$ is also constant,
and the $\xi_{i}$'s are all equal; but this is equivalent to
$\xi_{i}=0$. \bull

\section{2.}{Separation of variables for the $sl_{2}$ Gaudin-Calogero system}  

Let $X$ be the elliptic curve $\CC^{\times}/q^{\ZZ}$, with marked points
$z_{\a}q^{\ZZ}$, $\a=1,...,N$. 
In the $sl_{2}$ case, the Gaudin-Calogero system (which plays the role
of the Gaudin system in the present situation, cf. \refER, \refN) 
takes place in the
space 
$$
\eqalign{
\cM_{G}^{(0)}(X, z_{\a})=\{
(\cE_{(t,t^{-1})}, & {\rm\ parabolic\  structure\  at\ }z_{\a}{\rm\
given\ 
by\ }t_{\a}\in
\CC P^{1})\}/ \cr & [(t,t_{\a})\sim(qt,z_{\a}t_{\a}),
(t,t_{\a})\sim(t,ut_{\a}), u\in\CC^{\times}],
t\in\CC^{\times}] \cr
&
=\CC^{\times}\times (\CC P^{1})^{N}/\CC^{\times}\ltimes \ZZ^{N}.
}$$ 
Here $\cE_{(t_{1},\cdots,t_{n})}$ is the bundle on $X$
defined by
$\CC^{\times}\times\CC^{n}/[(z,\xi)\sim(qz,\diag(t_{i})\xi)]$, for
$t_{1},\cdots,t_{n}\in\CC^{\times}$. 

We consider then the space
$$
\eqalign{
\cM_{B}^{(0)}(X, z_{\a})=\{
 & (\cE_{(t,t^{-1})},  j:\cE_{t}\to \cE_{(t,t^{-1})},  {\rm \  par.\  str.\
given\  by\ } 
t_{\a}\in\CC P^{1}),t\in\CC^{\times}
\}/ \cr &[j\sim\la j,\la\in\CC^{\times}, 
(t,t_{\a})\sim(qt,z_{\a}t_{\a}),
(t,t_{\a})\sim(t,ut_{\a}), u\in\CC^{\times}];
}$$ 
it has a natural
projection $p$ to $\cM_{G}^{(0)}(X, z_{\a})$. Consider now the diagram 
$$
\matrix{
 & & & &{\cal Z}& & \cr
 & & &p_{1}\swarrow& &\searrow q_{1}& \cr
 & & & & & & \cr
 & & &\cM_{B}^{(0)}(X,z_{\alpha})& & &\{ {\cal E}\in {\rm
Pic}^{0}(X), \omega\in \PP\Hom (K_{X,z_{{\a}}}^{-1}
\otimes  {\cal E}^{-1},{\cal E})\}
\cr
 & & & & & & \cr
 & &p\swarrow & p_{0}\searrow& &\swarrow q_{0} &\searrow\pi \cr
 & & & & & & \cr
 & &\cM_{G}^{(0)}(X,z_{\alpha}) & &{\rm Pic}^{0}(X) & &
X^{(N)}\cr}
\leqno{(7)}
$$
where $K_{X,z_{\a}}=\Omega^{1}_{X}(\mathop{\sum}\limits_{\alpha
=1}^{N}({z}_{\alpha}))$, 
$p_{0}$ is the projection 
$$
{\rm class}(\cE_{(t,t^{-1})},j:\cE_{t}\to 
\cE_{(t,t^{-1})},{\rm \  par.\  str.})\mapsto\cE_{t^{-1}}
$$ 
($p_{0}$ associates
to $j$ its cokernel), $q_{0}$ associates $\cE$ to $(\cE,\omega)$, 
$\pi$ associates to $(\cE, \omega)$ the set of zeroes of $\omega$, and 
$\cZ$ is the incidence variety, defined by the conditions
that $\omega$ lifts to a morphism $j':
K^{-1}_{X,z_{\a}}\otimes\cE_{t^{-1}}\to
\cE_{(t,t^{-1})}$, compatible with the parabolic structure. Writing 
$$
\omega=\sum_{\a=1}^{N}u_{\a}{\theta(t^{-2}zz_{\a}^{-1})\over
{\theta(t^{-2})\theta(zz_{\a}^{-1})}}{dz\over z},
$$ 
the first component of $j'$ has to be
$\sum_{\a=1}^{N}u_{\a}t_{\a}{\dot\theta(zz_{\a}^{-1})\over \theta(zz_{\a}^{-1})
}{dz\over z}$, so that the incidence condition is
$\sum_{\a=1}^{N}u_{\a}t_{\a}=0$. 

(Fix our conventions for $\theta$- and $\wp$-functions:
$\theta(z)=\prod_{i\ge 0}(1-q^{i}z)\prod_{i>0}(1-q^{i}z^{-1})$, $\wp(\ln
z)=
-({\dot\theta\over\theta})^{\dot{}}(z)$, so $\wp(\tau)\sim
\tau^{-2}+\cdots$ for $\tau\to 0$; we denote $\dot
f(z)=z{{df}\over{dz}}$.)

The lift to $\CC^{\times}\times (\CC P^{1})^{N}$ of the Gaudin-Calogero
operators (\refER, \refN) are defined as follows: 
let $z$ be a formal variable, belonging to $X-\{z_{\a}\}$. We have 
$$
\eqalign{
L (z) & = e(z)f(z)+f(z)e(z)+{1\over 2}h(z)^{2}\cr
&= L_{0}+\sum_{\a=1}^{N}L_{\a}{\dot\theta\over\theta}(zz_{\a}^{-1})
+\sum_{\a=1}^{N}2\la_{\a}(\la_{\a}-1)\wp(\ln zz_{\a}^{-1})\cr &+
\sum_{\a=1}^{N}{1\over
2}h^{(\a)}(\sum_{\a=1}^{N}h^{(\a)})({\dot\theta\over\theta}(zz_{\a}^{-1}))^{2},
\cr}$$
where 
$e(z)=\sum_{\a=1}^{N}{\theta(t^{-2}zz_{\a}^{-1})\over
{\theta(t^{-2})\theta(zz_{\a}^{-1})}}e^{(\a)}$, 
$h(z)=2t^{2}{\pr\over{\pr t^{2}}}+2k{\dot\theta\over\theta}(t^{2})
+\sum_{\a=1}^{N}{\dot\theta\over\theta}
(zz_{\a}^{-1})h^{(\a)}$, 
$f(z)=\sum_{\a=1}^{N}{\theta(t^{2}zz_{\a}^{-1})\over
{\theta(t^{2})\theta(zz_{\a}^{-1})}}f^{(\a)}$, and
$e^{(\a)}=t_{\a}^{2}{\pr\over{\pr t_{\a}}}+2\la_{\a}t_{\a}$, 
$f^{(\a)}=-{\pr\over{\pr t_{\a}}}$, $h^{(\a)}=2(t_{\a}{\pr\over{\pr
t_{\a}}}+\la_{\a})$. 

Let us fix now complex numbers $\mu_{\a}$, $\a=0, \cdots, N$, with
$\sum_{\a=1}^{N}\mu_{\a}=0$ (this condition corresponds to the fact that
$\sum_{\alpha=1}^{N}L_{\alpha}$ belongs to the left ideal generated by
$\sum_{\alpha=1}^{N}h^{(\alpha)}$); consider on $\cM_{G}^{(0)}(X,
z_{\a})$, the
$\cD$-module (twisted by the quotient $\cL_{k,(\la_{\a})}$ of
$i^{-1}(\cL_{k}^{\boxtimes
2})\boxtimes\boxtimes_{\a=1}^{N}\cL_{\la_{\a}}$, $i:\Ker(s)\to X^{(2)}$,
$s:X^{(2)}\to X$ the sum mapping, $\cL_{k}$ a bundle of degree $k$ on $X$), 
$$
\cM_{(\mu_{\a})}=
\Big(\cD_{\cM_{G}^{(0)}(X, z_{\a})}\Big)_{\cL_{k,(\la_{\a})}}
/\sum_{\a=0}^{N}\Big(\cD_{\cM_{G}^{(0)}(X, z_{\a})}\Big)_{\cL_{k,(\la_{\a})}}
(L_{\a}-\mu_{\a}). 
$$ 
Consider then the operator on $X-\{z_{\a}\}$, 
$$
D_{(\la_{\a}),(\mu_{\a})}=2\Big(w{\pr\over{\pr w}}\Big)^{2}
-\mu_{0}-\sum_{\a=1}^{N}\mu_{\a}
{\dot\theta\over\theta}(wz_{\a}^{-1})+2\sum_{\a=1}^{N}\la_{\a}(\la_{\a}-1)
\wp(\ln wz_{\a}^{-1}),
$$ 
and the $\cD_{X}$-module
$$
\cE_{(\la_{\a}),(\mu_{\a})}=\cD_{X}/\cD_{X}D_{(\la_{\a}),(\mu_{\a})}.
$$
The fibration $\pi$ has fibers $\CC^{\times}$; we will twist inverse
images under $\pi$ by the function $C^{\sum_{\a=1}^{N}(\la_{\a}+1)}$
($C$ coordinate on the fiber). Also, we will work with $k=0$
(multiplication by $\theta(t^{2})^{k}$ taking us back to this case). 
We will show that:
 
\proclaim{Proposition} There is a homomorphism of $\cD$-modules
from the twisted inverse image 
$\pi^{*}\cE^{(N)}_{(\la_{\a}),(\mu_{\a})}$ to
$R(q_{1})_{*}p_{1}^{*}p^{*}
\cM_{(\mu_{\a})}[N]$, which is an isomorphism 
over $\pi^{-1}(X^{(N)}-\Delta)$
($\Delta$ is the diagonal part of $X^{(N)}$). 
\endproclaim

\demo Let 
$$
p_{1}: \CC^{\times}\times (\CC P^{1})^{N}\to
\CC^{\times}\times (\CC P^{1})^{N}/\CC^{\times}\ltimes\ZZ^{N}
$$ 
be the natural projection, then 
$$
p_{1}^{*}(\cM_{(\mu_{\a})})=\cD_{\CC^{\times}\times (\CC P^{1})^{N}}
/\sum_{\a}
\cD_{\CC^{\times}\times (\CC P^{1})^{N}}(L_{\a}-\mu_{\a})
+\cD_{\CC^{\times}\times (\CC P^{1})^{N}}(\sum_{\a}h^{(\a)}). 
$$

Because of the factor $\cD_{\CC^{\times}\times (\CC
P^{1})^{N}}(\sum_{\a}h^{(\a)})$, the $\cD$-module is constant along the
fibers of the action of $\CC^{\times}$. Its Radon transform 
is the $\cD$-module generated by $\bar
L_{\a}-\mu_{\a}$ and $\sum_{\a=1}^{N}\bar{h}^{(\a)}$, where 
$$
\eqalign{
 & \bar L_{0}  +\sum_{\a=1}^{N}\bar L_{\a}{\dot\theta\over\theta}(zz_{\a}^{-1})
+\sum_{\a=1}^{N}2\la_{\a}(\la_{\a}-1)\wp(\ln
zz_{\a}^{-1})+\sum_{\a=1}^{N}
{1\over 2}\bar h^{(\a)}(\sum_{\a=1}^{N}\bar
h^{(\a)})({\dot\theta\over\theta}(zz_{\a}^{-1}))^{2}\cr
&=\bar e(z)\bar f(z)+\bar f(z)\bar e(z)+{1\over 2}\bar h(z)^{2},
\cr}$$
$\bar e(z)=\sum_{\a=1}^{N}{\theta(t^{-2}zz_{\a}^{-1})\over
{\theta(t^{-2})\theta(zz_{\a}^{-1})}}\bar e^{(\a)}$, analogous formulae
for $\bar f(z)$ and $\bar h(z)$, with
$$
\bar e^{(\a)}=-[u_{\a}({\pr\over{\pr
u_{\a}}})^{2}+2(\la_{\a}+1){\pr\over{\pr u_{\a}}}], \bar
f^{(\a)}=u_{\a}, \bar h^{(\a)}=-2[u_{\a}{\pr\over{\pr
u_{\a}}}+(\la_{\a}+1)].$$ 

Consider the operator 
$$
\hat L(w_{i})=\bar L_{0}-\mu_{0}
+\sum_{\a=1}^{N}
{\dot\theta\over\theta}(w_{i}z_{\a}^{-1})(\bar L_{\a}-\mu_{\a})
+\sum_{\a=1}^{N}
({\dot\theta\over\theta}(w_{i}z_{\a}^{-1}))^{2} \lbreak
{1\over 2}\bar h^{(\a)}(\sum_{\a=1}^{N}\bar
h^{(\a)}).
$$ 
Let $\hat
e(w_{i})=\sum_{\a=1}^{N}{\theta(t^{-2}w_{i}z_{\a}^{-1})
\over{\theta(t^{-2})\theta(w_{i}z_{\a}^{-1})}}\bar e^{(\a)}$, etc. Let us
compute the difference
$$\eqalign{
-\hat L(w_{i})+(\hat e(w_{i})\hat f(w_{i})+\hat f(w_{i})\hat
e(w_{i})+{1\over 2}\hat
h(w_{i})^{2}-\mu_{0} & -\sum_{\a=1}^{N}\mu_{\a}{\dot\theta\over\theta}
(w_{i}z_{\a}^{-1})\cr &+2\la_{\a}(\la_{\a}-1)\wp(\ln
w_{i}z_{\a}^{-1})) 
;}
$$ it is equal to 
$$
\eqalign{
\sum_{\a,\b=1}^{N}{\theta(t^{-2}w_{i}z_{\a}^{-1})\over{\theta(t^{-2})
\theta(w_{i}z_{\a}^{-1})}}
 & [\bar e^{(\a)},
{\theta(t^{2}w_{i}z_{\b}^{-1})\over{\theta(t^{2})
\theta(w_{i}z_{\b}^{-1})}}]\bar f^{(\b)}
 \cr &+{1\over 2}
\sum_{\a,\b=1}^{N}{\dot\theta\over\theta}(w_{i}z_{\a}^{-1})
[\bar h^{(\a)}, {\dot\theta\over\theta}(w_{i}z_{\b}^{-1})]\bar h^{(\b)}. 
}\leqno{(8)}
$$
The first term of $(8)$ is the sum of $(9)$ and $(10)$, where
$$\eqalign{
(9)=-\sum_{\a,\b=1}^{N} & {\theta(t^{-2}w_{i}z_{\a}^{-1})\over{\theta(t^{-2})
\theta(w_{i}z_{\a}^{-1})}}
\big[2u_{\a}{\pr\over{\pr u_{\a}}}
({\theta(t^{2}w_{i}z_{\b}^{-1})\over
{\theta(t^{2})\theta(w_{i}z_{\b}^{-1})}})
{\pr\over{\pr u_{\a}}} \cr &+
u_{\a}({\pr\over{\pr u_{\a}}})^{2}
({\theta(t^{2}w_{i}z_{\b}^{-1})\over
{\theta(t^{2})\theta(w_{i}z_{\b}^{-1})}})\big]\cdot u_{\b}\cr
&= 
-\sum_{\a=1}^{N} {\theta(t^{-2}w_{i}z_{\a}^{-1})\over{\theta(t^{-2})
\theta(w_{i}z_{\a}^{-1})}}
2u_{\a}[\sum_{\b=1}^{N}{\pr\over{\pr u_{\a}}}
({\theta(t^{2}w_{i}z_{\b}^{-1})\over
{\theta(t^{2})\theta(w_{i}z_{\b}^{-1})}})u_{\b}
]{\pr\over{\pr u_{\a}}}\cr &
-
\sum_{\a=1}^{N}  {\theta(t^{-2}w_{i}z_{\a}^{-1})\over{\theta(t^{-2})
\theta(w_{i}z_{\a}^{-1})}}
2u_{\a}{\pr\over{\pr u_{\a}}}
({\theta(t^{2}w_{i}z_{\b}^{-1})\over
{\theta(t^{2})\theta(w_{i}z_{\b}^{-1})}})
{\pr\over{\pr u_{\a}}} \cr &
-
\sum_{\a=1}^{N} {\theta(t^{-2}w_{i}z_{\a}^{-1})\over{\theta(t^{-2})
\theta(w_{i}z_{\a}^{-1})}}
u_{\a}\sum_{\b=1}^{N}({\pr\over{\pr u_{\a}}})^{2}
({\theta(t^{2}w_{i}z_{\b}^{-1})\over
{\theta(t^{2})\theta(w_{i}z_{\b}^{-1})}})u_{\b}.
\cr}
$$
and 
$$
\eqalign{
(10)&=\sum_{\a,\b=1}^{N} {\theta(t^{-2}w_{i}z_{\a}^{-1})\over{\theta(t^{-2})
\theta(w_{i}z_{\a}^{-1})}} 
(-2)(\la_{\a}+1){\pr\over{\pr u_{\a}}}
 ({\theta(t^{2}w_{i}z_{\b}^{-1})\over{\theta(t^{2})
\theta(w_{i}z_{\b}^{-1})}})u_{\b}\cr
&
=\sum_{\a=1}^{N} 2(\la_{\a}+1)
{\theta(t^{-2}w_{i}z_{\a}^{-1})\over{\theta(t^{-2})
\theta(w_{i}z_{\a}^{-1})}} 
 {\theta(t^{2}w_{i}z_{\b}^{-1})\over{\theta(t^{2})
\theta(w_{i}z_{\b}^{-1})}}\cr
&
=\sum_{\a=1}^{N}2(\la_{\a}+1)[\wp(\ln t^{2})-\wp(\ln w_{i}z_{\a}^{-1})]
\cr}
$$

We have 
$$
\sum_{\b=1}^{N}{\pr\over{\pr u_{\a}}}
({\theta(t^{2}w_{i}z_{\b}^{-1})\over
{\theta(t^{2})\theta(w_{i}z_{\b}^{-1})}}
)u_{\b}=-{\theta(t^{2}w_{i}z_{\a}^{-1})\over
{\theta(t^{2})\theta(w_{i}z_{\a}^{-1})}},
$$
and
$$
\sum_{\b=1}^{N}({\pr\over{\pr u_{\a}}})^{2}
({\theta(t^{2}w_{i}z_{\b}^{-1})\over
{\theta(t^{2})\theta(w_{i}z_{\b}^{-1})}}
)u_{\b}=-2{\pr\over{\pr u_{\a}}}({\theta(t^{2}w_{i}z_{\a}^{-1})\over
{\theta(t^{2})\theta(w_{i}z_{\a}^{-1})}}).
$$
So 
$$
\eqalign{
(9)  =\sum_{\a=1}^{N}2u_{\a}{\theta(t^{2}w_{i}z_{\a}^{-1})\over{\theta(t^{2})
\theta(w_{i}z_{\a}^{-1})}} &
{\theta(t^{-2}w_{i}z_{\a}^{-1})\over{\theta(t^{-2})
\theta(w_{i}z_{\a}^{-1})}}  {\pr\over{\pr u_{\a}}}
\cr &=\sum_{\a=1}^{N}(-2)[\wp(\ln w_{i}z_{\a}^{-1})-\wp(\ln t^{2})
]u_{\a}{\pr\over{\pr u_{\a}}}.
}$$
The second term of $(8)$ is the sum of $(11)$ and $(12)$, with 
$$
(11)=2\sum_{\a,\b=1}^{N}{\dot\theta\over\theta}(w_{i}z_{\a}^{-1})
u_{\a}{\pr\over{\pr u_{\a}}}({\dot\theta\over\theta}(w_{i}z_{\b}^{-1}))
u_{\b}{\pr\over{\pr u_{\b}}}
$$
and
$$
\eqalign{
(12) 
&=
{1\over 2}\sum_{\a,\b=1}^{N}{\dot\theta\over\theta}(w_{i}z_{\a}^{-1})
2u_{\a}{\pr\over\pr u_{\a}}({\dot\theta\over\theta}(w_{i}z_{\b}^{-1}))
2(\la_{\b}+1)
\cr
&=
\sum_{\a=1}^{N}(\la_{\a}+1)\{-2w_{i}{\pr\over\pr w_{i}}
[{\dot\theta\over\theta}(w_{i}z_{\a}^{-1})]\}   
\cr
&=
\sum_{\a=1}^{N}2(\la_{\a}+1)\wp(\ln(w_{i}z_{\a}^{-1}))
\cr}$$

To compute $(11)$, we express the relation between the ${\pr\over{\pr
w_{i}}}$ and the ${\pr\over{\pr u_{\a}}}$: we have 
$$
\sum_{\a=1}^{N}u_{\a}{\theta(t^{2}zz_{\a}^{-1})\over{\theta(t^{2})
\theta(zz_{\a}^{-1})}}=
C{{\prod_{i=1}^{N}\theta(zw_{i}^{-1})}\over
{\prod_{\a=1}^{N}\theta(zz_{\a}^{-1})}},
$$
so 
$$\eqalign{
\sum_{\a=1}^{N}du_{\a}{\theta(t^{2}zz_{\a}^{-1})\over{\theta(t^{2})
\theta(zz_{\a}^{-1})}}
 & +{{dt}\over t}u_{\a}[{\dot\theta\over\theta}(t^{2}zz_{\a}^{-1})-
{\dot\theta\over\theta}(t^{2})]
{\theta(t^{2}zz_{\a}^{-1})\over{\theta(t^{2})
\theta(zz_{\a}^{-1})}} \cr
 & =
C{{\prod_{i=1}^{N}\theta(zw_{i}^{-1})}\over
{\prod_{\a=1}^{N}\theta(zz_{\a}^{-1})}}[{{dC}\over
C}-\sum_{i=1}^{N}{{dw_{i}}\over w_{i}}{\dot\theta\over\theta}
(zw_{i}^{-1})];\cr}
$$
by inspection of the pole at $z_{\a}$, and because of
$C{{\prod_{i=1}^{N}\theta(z_{\a}w_{i}^{-1})}\over
{\prod_{\b\ne\a}\theta(z_{\a}z_{\b}^{-1})}}=u_{\a}$, it follows that 
$$
du_{\a}=u_{\a}[{{dC}\over
C}-\sum_{i=1}^{N}{{dw_{i}}\over w_{i}}{\dot\theta\over\theta}
(z_{\a}w_{i}^{-1})],
$$
and so 
$$
w_{i}{\pr\over{\pr w_{i}}}
=-\sum_{\a=1}^{N}{\dot\theta\over\theta}
(z_{\a}w_{i}^{-1})
u_{\a}{\pr\over{\pr u_{\a}}}.
$$
Now 
$$
(11)=2\sum_{\b=1}^{N}-w_{i}{\pr\over{\pr w_{i}}}
({\dot\theta\over\theta}
(z_{\b}w_{i}^{-1}))u_{\b}{\pr\over{\pr u_{\b}}}
=2\sum_{\b=1}^{N}\wp(\ln w_{i}z_{\b}^{-1})u_{\b}{\pr\over{\pr u_{\b}}}. 
$$
Finally 
$$
(9)+(11)=2\wp(\ln t^{2})\sum_{\a=1}^{N}u_{\a}{\pr\over{\pr u_{\a}}}=
2\wp(\ln t^{2})C{\pr\over{\pr C}},
$$
and
$$
(10)+(12)=\sum_{\a=1}^{N}2(\la_{\a}+1)\wp(\ln t^{2}). 
$$

The sum of these terms is $2\wp(\ln t^{2})(C{\pr\over{\pr
C}}+\sum_{\a=1}^{N}(\la_{\a}+1))$;
 the term $C{\pr\over{\pr
C}}+\sum_{\a=1}^{N}(\la_{\a}+1)$ (equal to
$-{1\over 2}\sum_{\a=1}^{N}\bar{h}^{(\a)}$)
is set to zero in the twisted inverse image
$\pi^{*}\cE^{(N)}_{(\la_{\a}),(\mu_{\a})}$. 
It follows that 
$$
\eqalign{
\hat L(w_{i})=2[t^{2}{\pr\over{\pr t^{2}}}+
w_{i}{\pr\over{\pr
w_{i}}}+A(w_{i})
]^{2}
& -\mu_{0}-\sum_{\a=1}^{N}\mu_{\a}{\dot\theta\over\theta}(w_{i}z_{\a}^{-1})
\cr &
+\sum_{\a=1}^{N}2\la_{\a}(\la_{\a}-1)\wp(\ln w_{i}z_{\a}^{-1}),
}$$ 
since
$\hat e(w_{i})=0$, and $\hat h(w_{i})=2[t^{2}{\pr\over{\pr t^{2}}}+
k{\dot\theta\over\theta}(t^{2})+w_{i}{\pr\over{\pr
w_{i}}}+A(w_{i})]$, with 
$A(w_{i})=-\sum_{\a=1}^{N}(\la_{\a}+1){\dot\theta\over\theta}
(w_{i}z_{\a}^{-1})$. The addition of the term $A(w_{i})$  corresponds to
the twisting by a $GL(1)$-connection and does not change the
$PGL(2)$-oper, as in \refF. The addition of the term $t^{2} \pr/\pr (t^{2})$
corresponds to the fact that we are working here with the variables
$(w_{i},t)$, which are linked by the relation
$t^{2}\prod_{i=1}^{N}w_{i}=\prod_{\a=1}^{N}z_{\a}$ (mod. $q^{\ZZ}$). 
Our
statement follows as before. \bull

\section{Remark.}{} As in the rational case, one may remark that the
conditions on the $\mu_{\a}$'s to satisfy the Bethe ansatz equations, that
can be found in \refFV, can be translated into the condition on the
projective connection defined by them, to have a single-valued
solution $\psi(z)=
\prod\theta(za_{i})/\prod_{\a=1}^{n}\theta(zz_{\a}^{-1})^{\la_{\a}}$. 
Thanks to the Leray formulae for Radon transformation \refAS, one could
expect the Bethe eigenvectors to be expressed in the form 
$$
\eqalign{
\Psi(t_{1},\cdots,t_{N}) = &
\int_{\Gamma}{C^{-\sum_{\a=1}^{N}(\la_{\a}+1)}
\psi(w_{1})\cdots\psi(w_{N)}\over{(\sum_{\a=1}^{N}
u_{\a}t_{\a}})^{k}} \cdot \cr & \cdot 
\sum_{\a=1}^{N}(-1)^{\a}u_{\a}du_{1}\wedge
\cdots \wedge \check{du_{\a}}
\wedge
\cdots\wedge du_{N},
}$$
$1\le k-\sum_{\a=1}^{N}(\la_{\a}+1)\le N-1$, 
the integration being on a suitable cycle in $\CC P^{N}$. In the general
case, this formula should lead to the computation of the monodromy of
the Gaudin-Calogero system (by deformation of the cycle of
integration). It might be interesting to express this
monodromy representation directly in terms of the one of the projective
connection associated to the $\mu_{\a}$'s. 

\vskip 1truecm
\noindent
{\bf References}
\bigskip
\item{\refAS} A. d'Agnolo, P. Schapira, {\sl Quantification de Leray de
la dualit\'e projective,} C. R. Acad. Sci. Paris, t.139, s\'erie I (1994),
595-8. 

\item{\refBD} A.A. Beilinson, V.G. Drinfeld, {\sl Quantization of
Hitchin's fibration and Langlands program,} preprint. 

\item{\refD} V.G. Drinfeld, {\sl Two-dimensional representations of the
fundamental group of a curve over a finite field and automorphic forms
on $GL(2)$,} Amer. J. Math. 105 (1983), 85-114. 

\item{\refER} B. Enriquez, V. Rubtsov, {\sl  Hitchin systems, higher
Gaudin operators and $r$-matrices,} preprint alg-geom/9503010. 

\item{\refFG} F. Falceto, K. Gawedski, {\sl Unitarity of the
Knizhnik-Zamolodchikov-Bernard connection and the Bethe Ansatz for the
elliptic Hitchin system,} hep-th/9604094. 

\item{\refFFR} B.L. Feigin, E.V. Frenkel, N. Reshetikhin, {\sl Gaudin
model, Bethe ansatz and correlation functions at the critical level,}
Commun. Math. Phys. 166 (1), 27-62 (1995). 

\item{\refFV} G. Felder, A. Varchenko, {\sl  Integral representation of
solutions of the elliptic
Knizhnik--Zamolodchikov--Bernard equations,} preprint hep-th/9502165. 

\item{\refF} E.V. Frenkel, {\sl Affine algebras, Langlands duality and
Bethe ansatz,} Proc. ICMP-94, 606-42, International
Press (1995),
(q-alg/9506003). 

\item{\refN} N. Nekrasov, {\sl Holomorphic bundles and many-body systems,}
preprint hep-th/95-03157. 

\item{\refS} E.K. Sklyanin, {\sl Separation of variables in the Gaudin
model,} J. Sov. Math. 47 (1989), 2473-88.

\medskip
\medskip\medskip
\section{}{}

B.E.: Centre de Math\'{e}matiques, URA 169 
du CNRS, Ecole Polytechnique, 91128 Palaiseau, France

B.F.: Landau Inst. for Theor. Physics, Kosygina 2, GSP-1, 117940 Moscow
V-334, Russia

V.R.: ITEP, Bol. Cheremushkinskaya, 25, 117259, 
Moscow, Russia.
\bye